\documentclass[11pt]{article}

\usepackage{amsmath,amsfonts,amsthm,xypic,color,epsfig}
\usepackage{pst-all}
\usepackage{pstricks}
\usepackage[center,footnotesize]{caption}

\theoremstyle{plain}
\newtheorem{theorem}{Theorem}[section]
\newtheorem{proposition}[theorem]{Proposition}

\theoremstyle{definition}
\newtheorem{definition}[theorem]{Definition}

\newtheorem{lem}[theorem]{Lemma}

\newtheorem{cor}[theorem]{Corollary}
\newtheorem{ex}[theorem]{Example}

\newenvironment{renumerate}%
{%
\begin{enumerate}}%
{\end{enumerate}%
}%

\newenvironment{remark}%
{\vskip6pt%
\noindent%
{\it Remark.}}%
{\vskip6pt}

{\vskip6pt%
\noindent%
{\it Remarks}. %
\begin{renumerate}}%
{\end{renumerate}\vskip6pt}

\newcommand{\R}{\text{${\mathbb R}$}}
\newcommand{\C}{\text{$\mathbb C$}}

\newcommand{\Z}{\text{$\mathbb Z$}}

\newcommand{\X}{\text{$\mathbb X$}}

\newcommand{\J}{\text{$\mathcal{J}$}}
\newcommand{\M}{\text{$\mathcal{M}$}}
\newcommand{\N}{\text{$\mathcal{N}$}}

\newcommand{\OO}{\mathcal{O}}

\newcommand{\gO}{\text{$\Omega$}}

\newcommand{\del}{\text{$\partial$}}

\newcommand{\mc}[1]{\text{$\mathcal{#1}$}}

\newcommand{\into}{\longrightarrow}
\newcommand{\noqed}{\let\qed\relax}

\renewcommand{\Re}{\mathrm{Re}}
\renewcommand{\Im}{\mathrm{Im}}

\newcommand{\IP}[1]{\langle #1 \rangle}
\newcommand{\nhood}{neighbourhood}

\newcommand{\gcs}{generalized complex structure}

\newcommand{\gcss}{generalized complex structures}

\newcommand{\gcm}{generalized complex manifold}
\newcommand{\gcms}{generalized complex manifolds}

\newcommand{\wrt}{with respect to}
\renewcommand{\iff}{if and only if}

\newcommand{\eps}{\varepsilon}

\numberwithin{equation}{section}

\title{Blow-up of generalized complex 4-manifolds}

\author{Gil R. Cavalcanti\footnote{Utrecht University; g.r.cavalcanti@uu.nl.}\ \  and Marco Gualtieri\footnote{University of Toronto;  mgualt@math.utoronto.ca.}}

\date{}
\begin{document}

\definecolor{red}{rgb}{1,0,0}
\definecolor{green}{rgb}{0,1,0}
\definecolor{blue}{rgb}{0,0,1}
\definecolor{yellow}{rgb}{1,1,0}
\definecolor{pink}{rgb}{1,0,1}
\definecolor{cian}{rgb}{0,1,1}
\definecolor{darkblue}{rgb}{0,0.2,0.2}
\definecolor{brown}{rgb}{0.5,0.4,0.2}

\maketitle 

\abstract{We introduce blow-up and blow-down operations for
generalized complex 4-manifolds.  Combining these with a surgery
analogous to the logarithmic transform, we then construct
generalized complex structures on $n \C P^2 \# m \overline{\C P^2}$
for $n$ odd, a family of 4-manifolds which admit neither complex nor
symplectic structures unless $n=1$. We also extend the notion of a
symplectic elliptic Lefschetz fibration, so that it expresses a
generalized complex 4-manifold as a fibration over a two-dimensional
manifold with boundary.}

\section*{Introduction}

Generalized complex structures~\cite{Hi03,Gu03} are a simultaneous
generalization of complex and symplectic structures.  Since their
introduction, it has been natural to ask whether generalized complex
manifolds encompass a genuinely larger class than complex or
symplectic manifolds.  Indeed, the only obstruction for existence
known is that the underlying manifold must be almost
complex~\cite{Gu03}.  Generalized complex structures in dimension 2
are either complex or symplectic, so this question becomes
nontrivial first in real dimension 4.

In \cite{CG06}, the authors answered the above question in the
affirmative, by constructing a \gcs\ on $3 \C P^2 \# 19 \overline{\C
P^2}$.  This manifold does not have complex or symplectic
structures, due to Kodaira's classification of complex surfaces and
the fact that it has vanishing Seiberg--Witten invariants
\cite{Kod64,Ta94,Wi94}.

In this article, we develop blow-up and blow-down operations for
generalized complex manifolds.  We then show how this significantly
enlarges the list of manifolds which are generalized complex but not
complex or symplectic; in particular we prove that $n \C P^2 \#
m\overline{\C P^2}$ has a \gcs\ \iff\ it is almost complex, i.e. $n$
is odd.

The four-dimensional generalized complex manifolds which we consider
may be understood classically as symplectic structures which acquire
a singularity along a two-dimensional submanifold called the
\emph{complex} or \emph{type change} locus, which itself acquires a
complex structure. We prove that if a point is in the complex locus,
then it may be blown up, just as a point on a complex surface.  This
is done by proving a normal form theorem for neighbourhoods of such
points.

Just as in the complex case, blowing up a point introduces an
exceptional divisor, although in our case it is not a complex curve
but rather a two-dimensional submanifold which is Lagrangian away
from its intersection with the complex locus; such submanifolds are
called \emph{generalized complex branes}.  We then show that these
branes have standard tubular neighbourhoods, in analogy to
Weinstein's Lagrangian neighbourhood theorem.  This allows us to
show that spherical branes intersecting the complex locus at one
point may be blown down.

Finally, we apply our blow-up and blow-down operations to
generalized complex 4-manifolds obtained from symplectic fiber sums
of rational elliptic surfaces via the logarithmic transform
introduced in \cite{CG06}.  Describing the symplectic 4-manifolds as
Lefschetz fibrations, we use the method of vanishing cycles to find
spherical branes in the associated generalized complex
four-manifolds which may be blown down. Using Kirby calculus, the
resulting generalized complex manifolds are then shown to be
diffeomorphic to $n \C P^2 \# m\overline{\C P^2}$ for $n$ odd.

Motivated by these manipulations, we extend the notion of a (genus
1) symplectic Lefschetz fibration so that it applies to generalized
complex 4-manifolds.  The base for these generalized Lefschetz
fibrations may be any two-dimensional manifold with boundary.  The
one-dimensional boundary of the base corresponds precisely to the
complex locus in the 4-manifold. A more thorough study of these
fibrations is in progress.

This paper is organized as follows: in the first section, we
introduce generalized complex structures and prove a \nhood\ theorem
for a point in the complex locus; in Section \ref{sec:branes}, we
recall the definition of branes and prove a \nhood\ theorem for
generic branes in 4-manifolds; in Section \ref{sec:blow-up}, we show
that it is possible to blow up points in the complex locus as well
as blow down spherical branes intersecting the complex locus
transversally at a single point;  in Section \ref{sec:surgery}, we
recall the surgery introduced in \cite{CG06} and describe its effect
on Lefschetz fibrations; in the final section, we prove that $n \C
P^2 \# m \overline{\C P^2}$ has a \gcs\ if $n$ is odd.  This last
result requires a Kirby calculus computation, which we provide in
the Appendix.


\section{Generalized complex structures}\label{sec:gcss}
Given a closed 3-form $H$ on a manifold $M$, the {\it Courant
bracket}~\cite{Courant,SW01} of sections of  $TM\oplus T^*M$  is
defined by
\[[X+\xi,Y+\eta]_H = [X,Y] + \mc{L}_X \eta -\mc{L}_Y\xi -\frac{1}{2}d(\eta(X) - \xi(Y)) + i_Y i_X H.\]
The bundle $TM \oplus T^*M$ is also endowed with a natural symmetric
pairing of signature $(n,n)$:
\[ \IP{X+\xi,Y+\eta} = \tfrac{1}{2} (\eta(X) + \xi(Y)).\]
\begin{definition}
A {\it \gcs} on $(M,H)$ is a complex structure $\J$ on the bundle
$TM \oplus T^*M$ which preserves the natural pairing and whose
$+i$-eigenbundle is closed under the Courant bracket.
\end{definition}
Recall that the differential forms $\Omega^{\bullet}(M)$ carry a
natural spin representation for the metric bundle $TM \oplus T^*M$;
the Clifford action of $X+\xi \in TM \oplus T^*M$ on $ \rho \in
\Omega^{\bullet}(M)$ is
\[ (X+\xi) \cdot \rho  = i_X \rho + \xi \wedge \rho.\]
The $+i$-eigenbundle of $\J$ is a maximal isotropic subbundle of
$T_\C M\oplus T_\C^*M$, and we may use the correspondence between
maximal isotropics and pure spinors to encode $\J$ as a line
subbundle of the complex differential forms.

\begin{definition}
The {\it canonical bundle} of $\J$ is the complex line bundle
$K\subset\wedge^{\bullet}T_\C^*M$ annihilated by the
$+i$-eigenbundle of $\J$.
\end{definition}
As shown in~\cite{Gu03}, a complex differential form $\rho$ must
satisfy the following properties in order to be a local generator of
the canonical bundle of a generalized complex structure (and hence
determine $\J$ uniquely):
\begin{enumerate}
\item At each point, $\rho$ has the algebraic form
\begin{equation}\label{eq:pure}
\rho =  e^{B+ i \omega}\wedge \gO,
\end{equation}
where $B$ and $\omega$ are real 2-forms and $\gO$ is a decomposable
complex form.
\item At each point, $\rho$ satisfies the nondegeneracy condition
\begin{equation}\label{eq:LcapLbar}
(\rho,\overline{\rho}) = \gO \wedge \overline{\gO} \wedge
(2i\omega)^{n-k} \neq 0.
\end{equation}
\item The form $\rho$ is integrable, in the sense
\[
d\rho + H\wedge\rho = (X+\xi)\cdot\rho,
\]
for some section $X+\xi$ of $TM\oplus T^*M$.
\end{enumerate}
The first condition is equivalent to the fact that $\rho$ must be a
pure spinor, which is a pointwise algebraic condition.  The second
condition derives from the transversality of the $\pm
i$-eigenbundles of $\J$, and involves the natural Spin-invariant
pairing of differential forms with values in the top degree forms:
\[
(\rho,\sigma)= [\rho^\top\wedge\sigma]_{top}.
\]
(Here, $\rho^\top$ denotes the reversal anti-automorphism of forms).
We see from this condition that the volume form
$i^{-n}(\rho,\overline\rho)$ defines a canonical global orientation
on any $2n$-dimensional generalized complex manifold.  We also
derive the fact that at each point of a generalized complex
manifold, $\ker \gO\wedge\overline{\gO}$ is a subspace of the real
tangent space with induced symplectic structure and transverse
complex structure.

%

\begin{definition}
Let $\J$ be a \gcs\ and $e^{B+i \omega} \wedge \gO$ a generator of
its canonical bundle at a point $p$. The {\it type} of $\J$ at $p$
is the degree of $\gO$ and the {\it parity} of $\J$ is the parity of
its type.
\end{definition}
We will shortly see examples where the type of a \gcs\ jumps along
loci in the manifold.  However, its parity must clearly remain
constant on connected components of $M$.

\begin{ex}\label{ex:complex} Let $(M^{2n},I)$ be a complex manifold. Then
the following operator on $TM \oplus T^*M$ is a generalized complex
structure:
\[ \J_I  = \begin{pmatrix} -I & 0 \\ 0 & I^* \end{pmatrix}\]
The $+i$-eigenspace of $\J_I$ is $T^{0,1}M \oplus T^{*{1,0}}M$,
which annihilates the canonical bundle $K=\wedge^{n,0}T^*M$ and is
therefore of type $n$. The orientation induced by the \gcs\ is the
same as the one induced by the underlying complex structure.
\end{ex}

\begin{ex}\label{ex:symplectic}
Let $(M,\omega)$ be a symplectic manifold. Then
\[\J_{\omega} = \begin{pmatrix}0& -\omega^{-1} \\ \omega & 0 \end{pmatrix}\]
is a generalized complex structure with $+i$-eigenspace $\{X -
i\omega(X):X \in T_{\C}M\}$ and canonical bundle generated by the
differential form $e^{i \omega}$.  Symplectic structures, therefore,
have type zero. The orientation induced by this \gcs\ is the same as
the orientation induced by the symplectic structure.
\end{ex}

\begin{ex}\label{typechang}
Let $(z,w)$ be standard complex coordinates on $\C^2$. The complex
differential form
\[
\rho = w + dw\wedge dz
\]
may be expressed as $\rho = w\exp(w^{-1}dw\wedge dz)$ when $w\neq 0$
and as $\rho = dw\wedge dz$ when $w=0$. Hence it is a pure spinor
by~\eqref{eq:pure}. Further, it is nondegenerate since
$(\rho,\overline\rho)=dw\wedge d\bar w\wedge dz\wedge d\bar z\neq
0$. Finally, we see that $\rho$ is integrable, since
\[
d\rho = -\del_z\cdot\rho.
\]
Hence $\rho$ defines a generalized complex structure on $\C^2$ which
undergoes type change: it has complex type (type 2) along the locus
$w=0$, and symplectic type (type 0) elsewhere.
\end{ex}

\begin{ex}\label{ex:B-field}
Any real 2-form $B$ gives rise to an orthogonal transformation of
$TM \oplus T^*M$ via $e^B:X+\xi\mapsto X+\xi-i_X B$.  This
transformation induces an isomophism between the $H$-Courant bracket
and the $H+dB$-Courant bracket, hence it acts by conjugation on any
given \gcs\ $\J$ on $(M,H)$, producing a new one $e^{-B}\J e^B$ on
$(M,H+dB)$. The induced action on the canonical bundle is simply
\[
K\mapsto e^B\wedge K= (1 + B + \tfrac{1}{2}B\wedge B + \cdots)\wedge
K.
\]
\end{ex}
Example~\ref{ex:B-field} indicates that there are symmetries of the
Courant bracket beyond the usual diffeomorphisms.  We now use this
to define morphisms between generalized complex manifolds.

\begin{definition}
Let $\M_i=(M_i,H_i)_{i=1,2}$ be manifolds equipped with closed
3-forms. Then $\Phi = (\varphi,B)\in C^\infty(M_1,M_2)\times
\Omega^2(M_1,\R)$ is called a generalized map $\M_1\rightarrow\M_2$
when $\varphi^*H_2- H_1 = dB$. When $\varphi$ is a diffeomorphism we
call $\Phi$ a {\it $B$-diffeomorphism}.
\end{definition}
A generalized map establishes a correspondence between the tangent
and cotangent bundles which is neither covariant nor contravariant.
We say that $X + \xi\in TM_1\oplus T^*M_1$ is $\Phi$-related to
$Y+\eta\in TM_2\oplus T^*M_2$, and write
\[
X+\xi \sim_{\Phi}Y+\eta,
\]
when $Y = \varphi_* X$ and $\xi = \varphi^*\eta + i_XB$.

\begin{definition}
Let $\Phi:\M_1\rightarrow \M_2$ be a generalized map, and let $\J_i$
be generalized complex structures on $\M_i$. Then $\Phi$ is
{\it holomorphic}\footnote{This notion of morphism essentially coincides
with that described in~\cite{Crainic} and~\cite{GuBi}.} when
\[
\J_1(X+\xi) \sim_{\Phi}\J_2(Y+\eta)
\]
for all $(X+\xi)\sim_{\Phi} (Y+\eta)$. When $\Phi$ is a
$B$-diffeomorphism, we say $\J_1,\J_2$ are {\it isomorphic}.
\end{definition}
This notion of morphism specializes to a holomorphic map if the
$\J_i$ are usual complex structures, and to a symplectomorphism if
the structures are symplectic. In the case of an isomorphism, we see
directly that $\J_1 = e^{-B}(\varphi^*\J_2) e^{B}$, while on
canonical bundles the isomorphism yields
\[
K_1 = e^B\varphi^*K_2.
\]

We now construct further examples of \gcss\ by deforming the usual
complex manifolds from Example~\ref{ex:complex}.
\begin{theorem}[Gualtieri \cite{Gu03}]\label{theo:deformations}
Any holomorphic Poisson bivector $\beta$ on a complex manifold
$(M,I)$ deforms the complex structure into a generalized complex
structure $\J_\beta$, with canonical bundle
\begin{equation}\label{kbeta}
K_\beta = e^\beta\Omega^{n,0},
\end{equation}
where $\beta$ acts by interior product.
\end{theorem}
The action of $\beta$ on $TM\oplus T^*M$ giving rise
to~\eqref{kbeta} is  $e^P: X+\xi\mapsto X -P(\xi) + \xi$, for
$P=\beta+\overline\beta$. Hence the deformed complex structure on
$TM\oplus T^*M$ is
\begin{equation}\label{pandj}
 \J_\beta =e^{P}\J_Ie^{-P} = \begin{pmatrix} -I & Q \\ 0 & I^* \end{pmatrix},
\end{equation}
where $Q=-4\Im(\beta)$.

Generalized complex structures obtained by deformation in this way
do not necessarily have constant type over $M$: the deformed
structure has type equal to the corank of $\beta$, which may vary
along the manifold. We now investigate several examples of deformed
complex surfaces, where the resulting generalized complex structure
has generic type zero (symplectic type) jumping to type $2$ (complex
type) along the vanishing locus of $\beta$, an anticanonical
divisor.

\begin{ex}\label{ex:linb}
Let $\Sigma$ be a Riemann surface and $\pi:L\rightarrow \Sigma$ be a
holomorphic line bundle.  The total space of $L$ is a complex
surface $S$, and hence any holomorphic bivector field on $S$ is
automatically Poisson.  Using the fact that $TS$ is an extension of
$\pi^*T\Sigma$ by $\pi^*L$, we see that $\wedge^2 TS =
\pi^*T\Sigma\otimes\pi^*L$.  By a fibrewise Taylor expansion about
the zero section, we obtain a filtration of $H^0(S, \wedge^2 TS)$ by
sections of polynomial degree at most $k$ along the fibers:
\[
H^0_{k}(S,\wedge^2 TS) = \bigoplus_{i=0}^k H^0(\Sigma,
L^{1-i}\otimes T\Sigma).
\]
Taking $i=1$ above, we obtain a holomorphic bivector vanishing to
order 1 along the zero section as long as $T\Sigma$ is trivial, i.e.
$\Sigma$ is an elliptic curve. Hence we obtain a generalized complex
structure of generic type $0$, jumping to complex type along an
elliptic curve, irrespective of the line bundle $L$.

Taking $i=2$ above, and setting $L = T\Sigma$, we obtain a
holomorphic bivector vanishing to order 2 along the zero section,
giving rise to a generalized complex structure which undergoes type
change along the zero section of $T\Sigma$, irrespective of
$\Sigma$.
\end{ex}

The preceding example suggests that the nature of the type change
locus depends on a certain order of vanishing. Indeed, a generalized
complex structure of generic type 0 will undergo type change
precisely where the projection of $K\subset \wedge^\bullet T^*_\C M$
to $\wedge^0 T^*_\C M = \C$ vanishes.  In other words, this
projection defines a section $s\in C^\infty(K^*)$ and the type
change occurs at the zero locus of $s$.  For a 4-dimensional
manifold, this is the only possible type change, since type 2 is
maximal.
\begin{definition}
A point $p$ in the complex (type 2) locus of a generalized complex
4-manifold is called {\it nondegenerate} if it is a nondegenerate
zero of the section $s\in \Gamma(K^*)$.  If $p$ is a zero of order
$a$ of the section $s$, then we call $p$ a degenerate complex point
of order $a$.
\end{definition}
In the remainder of this section, we show that
Example~\ref{typechang} provides a normal form for a neighbourhood
of any nondegenerate complex point.  We may view this model as a
deformation of the standard $\C^2$ by the bivector $\beta = w
\del_w\wedge\del_z$, since
\begin{equation}\label{eq: point local form}
 \rho =e^\beta\cdot dw \wedge dz =  w + dw\wedge dz.
\end{equation}

\begin{theorem}\label{prop:local form}
Let $(M,\J)$ be a generalized complex 4-manifold and let $p \in M$
be a nondegenerate complex point. Then $p$ has a neighbourhood which
is $B$-diffeomorphic to a \nhood\ of the origin in $\C^2$ with the
\gcs\ determined by \[\rho=w + dw \wedge dz.\]
\end{theorem}
\begin{proof}
Let $\rho= \rho_0 + \rho_2 + \rho_4$, $deg(\rho_i) = i$, be a
trivialization of $K$ in a neighbourhood of $p$ with $\rho_0(p)=0$.
Since $\rho$ is annihilated by the $+i$-eigenbundle of $\J$, there
is a unique real section $X+\xi\in C^\infty(TM\oplus T^*M)$ such
that
\begin{equation}\label{intpr}
d\rho = (X+\xi)\cdot\rho.
\end{equation}
Nondegeneracy implies that $d\rho_0|_p \neq 0$.  Hence
condition~\eqref{intpr} implies that $X(p) \neq 0$, and therefore
$X$ is nonzero in a \nhood\ of $p$. So we can parametrize a \nhood\
of $p$ by $(v,x) \in \R^3 \times \R$, with $X = \del_x$. Then the
closed 2-form
\[ B(v,x) = \int_0^x d\xi(v,t)dt\]
is such that
\[ d(i_XB - \xi )= \mc{L}_X B - d\xi = 0.\]
Therefore, $i_XB - \xi = df$ for some function $f$.  Finally we
obtain that
\[d(e^{f+B}\rho) = (df +X +\xi -i_XB)\cdot e^{f+B}\rho = X\cdot e^{f+B}\rho,\]
so that $d\tilde\rho = X\cdot\tilde\rho$ for $\tilde\rho =
e^{f+B}\rho$.  Therefore, we see that $\J$ is $B$-diffeomorphic,
near $p$, to a generalized complex structure, which we henceforth
denote $\J$, whose canonical bundle is generated by a form $\rho$
which satisfies $d\rho = X\cdot\rho$ for a real vector field $X$. An
immediate consequence of this is that $\mc{L}_X\rho = 0$, and hence
$\J$ is invariant in the $X$ direction.

Now consider the real section $\J X = Y+\eta$: since $\J$ is
$X$-invariant, we see that $\mc{L}_X(Y+\eta)=0$, implying that
$[X,Y]=0$ and $X\cdot d\eta = 0$ (since $\eta(X) = \IP{X,\J X}=0$,
by orthogonality of $\J$).

Since $\J$ is a complex structure at $p$, the real vector field $Y$
is nonvanishing near $p$. Since $[X,Y]=0$, we may choose coordinates
$(v,x,y) \in \R^2\times \R \times \R$ near $p$ such that $X =
\del_x$ and $Y = \del_y$.  Then define the closed 2-form
\[B = \tilde B + dy\wedge (\eta - i_Y\tilde B),\]
where $\tilde B$ is the closed 2-form defined by
\[\tilde B(v,x,y) = \int_0^y d\eta(v,x,t)dt.\]
The form $B$ is constructed precisely so that $i_X B = 0$ while $i_Y
B  = \eta$.  Therefore we obtain
\begin{align*}
e^{-B}\J e^{B} X = Y,
\end{align*}
showing that $\J$ is $B$-diffeomorphic, near $p$, to a generalized
complex structure (henceforth denoted $\J$) with generator $\rho$
satisfying $d\rho = X\cdot\rho$ and such that $\J X = Y$, for
nonvanishing real vector fields $X,Y$.

A direct result is that $X-iY$ lies in the $+i$-eigenbundle of $\J$,
which annihilates $\rho$.  In particular, $(X-iY)\cdot\rho_4=0$,
which implies $\rho_4=0$.  Therefore, we have $\rho = \rho_0 +
\rho_2$, with nondegeneracy guaranteeing
$\rho_2\wedge\overline\rho_2\neq 0$ and integrability giving
$d\rho_2=0$. Therefore $\rho_2$ determines a complex structure in
the neighbourhood of $p$. Integrability also implies
$d\rho_0\wedge\rho_2=0$, meaning $\rho_0$ is a holomorphic function;
define $w=\rho_0$ and choose $z$ so that $\rho_2=dw\wedge dz$. These
coordinates therefore render $\J$ into the desired normal form.
\end{proof}

Theorem~\ref{prop:local form} allows us to determine what structures
are inherited by the complex locus from the ambient generalized
complex geometry.  Let $U_i$ be neighbourhoods for which the above
theorem holds, with generalized complex structures defined by
$\rho_i=w_i+dw_i\wedge dz_i$.   On the overlaps $U_i\cap U_j$ we
have
\begin{equation}\label{locrho}
\rho_i = g_{ij}e^{B_{ij}}\rho_j,
\end{equation}
for $g_{ij}:U_i\cap U_j\rightarrow \C^*$ and
$B_{ij}\in\Omega^2_{cl}(U_i\cap U_j,\R)$ smooth \v{C}ech cocycles.
In degrees $0$ and $2$, this equation becomes
\begin{align}
w_i &= g_{ij}w_j\label{loc1}\\
dw_i\wedge dz_i &= g_{ij} dw_j\wedge dz_j +
w_jg_{ij}B_{ij}.\label{loc2}
\end{align}
Differentiating~\eqref{loc1} and subtracting from~\eqref{loc2}, we
obtain
\begin{equation}\label{transition}
g_{ij} dw_j\wedge ( dz_j - dz_i) = w_j(g_{ij}B_{ij} + dg_{ij}\wedge
dz_i),
\end{equation}
which vanishes on the complex locus $\Sigma$, defined by $w_j=0$.
Expanding~\eqref{transition} along $\Sigma$, we obtain
\begin{equation}\label{vann}
g_{ij} dw_j\wedge ( dz_j - \tfrac{\del z_i}{\del z_j}dz_j
-\tfrac{\del z_i}{\del \bar z_j}d\bar z_j  - \tfrac{\del z_i}{\del
\bar w_j}d\bar w_j )\big|_\Sigma = 0,
\end{equation}

This implies that $z_j$ is a holomorphic function of $z_i$ on
$\Sigma$, and furthermore $\del z_j/\del z_i = 1$. Hence the complex
locus, where it is nondegenerate, inherits a complex structure with
a distinguished trivialization of its holomorphic tangent bundle.
This trivialization $\del_z$ corresponds precisely to the vector
field $X+iY$ in the proof of Theorem~\ref{prop:local form}.


Taking the Lie derivative of~\eqref{loc2} in the $\overline z_j$
direction and restricting to the complex locus, we see that $\del
g_{ij}/\del \overline z_j=0$, showing that the conormal bundle of $\Sigma$ inherits a holomorphic structure.
Summarizing, we obtain the following result, extending work
in~\cite{CG06}.

\begin{cor}\label{cystr}
Let $\Sigma$ be the set of nondegenerate type changing points of a
4-dimensional generalized complex manifold $M$. Then $\Sigma$ is a
smooth 2-dimensional submanifold and inherits a holomorphic
structure (i.e. it is a Riemann surface) as well as a distinguished
trivialization $Z$ of its holomorphic tangent bundle (of course,
this determines a holomorphic differential $\Omega = Z^{-1}$). It
follows immediately that any compact component of $\Sigma$ must be
an elliptic curve. Furthermore, the conormal
bundle $N^*\Sigma$ inherits the structure of a
holomorphic bundle.
\end{cor}

\begin{remark}
The holomorphic vector field $Z = X+iY$ induced on the nondegenerate
complex locus $\Sigma$ has no canonical extension to the whole
4-manifold $M$. However, there is a natural extension of $Y$ to a
global class in the Poisson cohomology of $M$ with respect to a real
Poisson structure $P$ obtained from the generalized complex
structure (see~\cite{Gu7} for details).  In fact this extension is
the \emph{modular class} of the Poisson structure $P$ in the sense
of Weinstein~\cite{Weinmod}.
\end{remark}


\section{Branes}\label{sec:branes}

Generalized complex branes~\cite{Gu7} are a natural class of
submanifolds defined for any generalized complex manifold. For usual
complex manifolds the definition specializes to the notion of
complex submanifold.  In the symplectic case, Lagrangian
submanifolds provide examples.  It is appropriate to call these
\emph{(D)-branes} since they provide boundary conditions in
topological open string theory.

\begin{definition}
A {\it brane} in a \gcm\ $(M,H,\J)$ is a submanifold $\iota:\Sigma
\hookrightarrow M$ together with a 2-form $F \in \Omega^2(\Sigma)$
satisfying $dF = \iota^*H$ and such that the subbundle
$\tau_F\subset (TM\oplus T^*M)|_N$, defined by
\[\tau_F = \{X + \xi \in T\Sigma \oplus T^*M: \iota^*\xi = i_XF\},\]
is invariant under $\J$ (i.e. $\tau_F$ is a complex subbundle).
\end{definition}
Note that $\tau_F$ is a maximal isotropic subbundle isomorphic to
$T\Sigma\oplus N^*\Sigma\subset TM\oplus T^*M$.  When $(M,\J)$ is a
$2n$-manifold of symplectic type, with canonical bundle generated by
$e^{B+i\omega}$, the branes of lowest dimension are Lagrangian
submanifolds $\iota: \Sigma\rightarrow M$ with $F=\iota^*B$. As
shown in~\cite{Gu03}, there may exist non-Lagrangian branes on a
symplectic $2n$-manifold, of dimension $n + 2k$, corresponding to
the coisotropic A-branes discovered by Kapustin and
Orlov~\cite{Kap}. On the other hand, when $(M,\J)$ is a usual
complex structure, then $(\Sigma,F)$ is a brane if and only if
$\Sigma$ is a complex submanifold and $F$ is a $(1,1)$-form.

It follows from these extremal cases that even generalized complex
4-manifolds may only have $0$, $2$, or $4$-dimensional branes.
Branes of dimension zero coincide with points in the complex locus,
since points of a symplectic manifold are not branes.   It is shown
in~\cite{Gu7} that branes of dimension $4$ only occur when $(M,\J)$
is a deformation of a complex manifold via
Theorem~\ref{theo:deformations}.  We therefore concentrate on the
description of branes of dimension 2.

\begin{proposition}\label{prop:lagrangian brane}
Let (M,H,\J) be an even generalized complex 4-manifold, and let
$\Sigma \subset M$ be a 2-dimensional submanifold intersecting the
complex locus transversally at nondegenerate points. Then there
exists $F\in\Omega^2(\Sigma)$ such that $(\Sigma,F)$ is a brane
\iff\ $\Sigma$ is Lagrangian away from the complex locus.
Furthermore, $F$ is unique when it exists.
\end{proposition}
\begin{proof}\label{prop:lagrangian}
One implication follows from the description of branes for complex
and symplectic structures.  For the other implication, let $U$ be
the dense open set where the structure is of symplectic type, i.e.
given by the form $e^{B+i \omega}$. Then $F = B|_{\Sigma}$ is the
unique $2$-form on $\Sigma \cap U$ such that $\tau_F$ is
$\J$-invariant over $\Sigma\cap U$. Since $\J$-invariance is a
closed condition, $\tau_F$ will be $\J$-invariant over all of
$\Sigma$ as long as $F$ extends smoothly to $\Sigma$. Therefore it
remains to show that $F$ extends smoothly to the points where
$\Sigma$ intersects the complex locus.

Near a nondegenerate complex point, Proposition \ref{prop:local
form} provides the following normal form for $\J$: it is defined by
\[\rho= w +dw \wedge dz,\]
with $w=0$ defining the complex locus.  In coordinates $(w,z) =
(x+iy,u+iv)$, we have that $\rho = we^{B+i\omega}$, for
\begin{align*}
B& = \tfrac{1}{x^2+y^2}(x(dx\wedge du- dy\wedge dv )+ y( dx\wedge dv + dy \wedge du)),\\
\omega &= \tfrac{1}{x^2+y^2}(x(dx\wedge dv + dy\wedge du) - y(
dx\wedge du - dy\wedge dv)).
\end{align*}
Since $\Sigma$ is transversal to the complex locus, we can
parametrize it as $\X(x,y)=(x,y,u(x,y),v(x,y))$ where $u$ and $v$
vanish at 0. Since $\Sigma$ is Lagrangian, we have
$\omega(\X_x,\X_y)=0$, yielding
\begin{equation}\label{eq:lagrangian}
x(v_y -u_x) - y(u_y + v_x) = 0.
\end{equation}
In particular, there exists a smooth function $f:\R^2 \into \R$ such
that $u_y + v_x = x \,f$ and we can compute
\begin{align*}
B(\X_x,\X_y)&= \tfrac{1}{x^2+y^2}(x(u_y+v_x) + y(v_y - u_x))\\
&=\tfrac{1}{x^2 + y^2}(x(u_y+v_x) +\tfrac{y^2}{x}(u_y+v_x))\\
&=\tfrac{1}{x}(u_y+v_x) = f,
\end{align*}
where in the second equality we used~\eqref{eq:lagrangian}. Hence the
restriction of $B$ to $\Sigma$ extends smoothly to $w=0$, showing
$F$ is well defined on all of $\Sigma$.
\end{proof}

In the case that $\J$ is obtained by deforming a complex structure
via Theorem~\ref{theo:deformations}, there is a rich source of
examples of such $2$-dimensional branes; in fact, curves in the
original complex surface remain branes after the deformation, as we
now show.
\begin{proposition}\label{prop:branes}
Let $\J_\beta$ be the generalized complex structure obtained by
deforming the complex surface $(M,I)$ by the holomorphic bivector
$\beta$. Then any smooth complex curve $\Sigma \subset M$ with
respect to $I$ is a brane for $\J_\beta$, with $F=0$.
\end{proposition}
\begin{proof}
$\Sigma$ is a curve for $I$, hence
\[\tau_0= \{X + \xi \in T\Sigma \oplus T^*M: \xi|_\Sigma = 0\} = T\Sigma \oplus N^*\Sigma\]
is invariant under $\J_I$. By Theorem~\ref{theo:deformations}, we
have $\J_\beta = e^P\J_Ie^{-P}$, for $P=\beta+\overline\beta$. The
result then follows from the fact that $e^P \tau_0 =\tau_0$, which
we now show. If $X + \xi \in \tau_0$, then $e^P(X + \xi) = X -P(\xi)
+\xi $.  But $P(N^*\Sigma)\subset T\Sigma$, i.e. $\Sigma$ is
$P$-coisotropic, since $\Sigma$ is of type $(1,1)$ with respect to
$I$ while $P$ is of type $(2,0)+(0,2)$.
\end{proof}

Example~\ref{ex:linb} provides examples of $\beta$-deformed complex
structures on the total space of a holomorphic line bundle $\pi:
L\rightarrow \Sigma$ over a Riemann surface; holomorphic bivectors
of homogeneous degree $i$ along the fibres are again given by
\begin{equation}\label{ordvan}
H^0(\Sigma, K^{-1}_\Sigma\otimes L^{1-i}).
\end{equation}
Before deformation, complex curves in these examples include the
zero section $\Sigma$, as well as any fiber $\pi^{-1}(p),\
p\in\Sigma$. Therefore, by Proposition~\ref{prop:branes}, all these
give examples of $2$-branes for $\J_\beta$. Taking $i=0$, we obtain
examples where the zero section $\Sigma$ is a generically Lagrangian
brane, as follows.

\begin{ex}\label{ex:line bundles}
Let $\Sigma$ be a Riemann surface, and $D = \sum a_i p_i,\ a_i>0$ an
effective divisor on $\Sigma$, for $p_1,\ldots, p_n$ points in
$\Sigma$. Let $\OO(D)$ be the associated holomorphic line bundle,
with section $\beta\in H^0(\Sigma,\OO(D))$ such that $D=(\beta)$.
By~\eqref{ordvan}, $\beta$ defines a deformation of the complex
structure on the total space of $L = K_{\Sigma}(D)$ which is
constant along the fibres and vanishes precisely on the fibres
$\pi^{-1}(p_i)$ above the divisor. Hence the zero section $\Sigma$
is a brane for $\J_\beta$ which is generically Lagrangian but
intersects the complex locus transversally at the points
$p_i\in\Sigma$; at these points the complex locus is degenerate of
order $a_i$.

This structure is easily described in coordinates: let
$U_0=\Sigma\backslash\{p_1,\ldots,p_n\}$ and let
$z_i:U_i\rightarrow\C$ be coordinates such that $z_i(p_i)=0$ and
$U_i\cap U_j=\emptyset$ for nonzero $i,j$ unless $i=j$. The bundle
$\OO(D)$ is taken to be trivial on the $U_i$, with transition
functions $f_{i0}(z_i) = z_i^{a_i}$.  The defining section for the
divisor $D$ is given by the functions $z_i^{a_i}$ (in $U_i,\ i\neq
0$) and the constant function $1$ in $U_0$.

The line bundle $L=K_\Sigma(D)$ is described by tensoring the above
trivialization of $\OO(D)$ with $K_\Sigma=T^*\Sigma$. To specify the
generalized complex structure, we give differential forms $\rho_i$
on $\tilde U_i = T^*U_i$ and cocycles $g_{ij},B_{ij}$ such
that~\eqref{locrho} holds.  Let $w_i$ be the canonically conjugate
coordinate to $z_i$ on $T^*U_i$, let $\Omega=B+i\omega$ be the
natural holomorphic symplectic form on $T^*\Sigma$, i.e.
$\Omega|_{U_i} = dz_i\wedge dw_i$, and define
\begin{align*}
\rho_0 &= e^{i\omega}|_{\tilde U_0},\\
\rho_i &= z_i^{a_i} + dz_i\wedge dw_i,\\
g_{i0}(z_i) &= z_i^{a_i},\\
B_{i0}&=B|_{\tilde U_{i}\cap \tilde U_0}.
\end{align*}
Note that we glue $T^*U_i$ to $T^*U_0$ via the diffeomorphism
$\varphi_{i0}:(z_i,w_i)\mapsto (z_i,z_i^{a_i}w_i)$, due to the
twisting by $D$. Therefore we have $\rho_i =
\varphi_{i0}^*z_i^{a_i}(1+dz_i\wedge dw_i)=
\varphi_{i0}^*g_{i0}e^{B_{i0}}\rho_0$ in $U_0\cap U_i$, as required.
\end{ex}

The generalized complex structure in the previous example does not
depend on the holomorphic structure of the initial Riemann surface
$\Sigma$, or on the locations of the points $p_i$, since a
(orientation-preserving) diffeomorphism $\psi:\Sigma\rightarrow
\Sigma'$ may always be chosen to send $p_i$ to $p_i'$ and to take
$z_i$ to $z_i'$ in a neighbourhood of the points $p_i, p_i'$. The
symplectic form on $T^*\Sigma$ is diffeomorphism invariant, hence we
obtain $\psi^*\rho_i' = \rho_i$ for all $i$, yielding an isomorphism
of generalized complex structures.  This suggests the following more
general example.

\begin{ex}\label{normalforms}
Let $\Sigma$ be a real smooth 2-manifold and $D = \sum a_i p_i,\
a_i>0,$ be a smooth effective divisor, i.e. a positive integer
linear combination of points $p_i$ of $\Sigma$.  Choose complex
coordinates $z_i$ in neighbourhoods $U_i\subset \Sigma$ of $p_i$ and
canonically conjugate coordinates $w_i$ in $\tilde U_i = T^*U_i$, so
that $dz_i\wedge dw_i = B_i + i\omega$, where $\omega$ is the
canonical real symplectic form on $T^*\Sigma$.  Then setting $\tilde
U_0=T^*(\Sigma\backslash\{p_1,\ldots,p_n\})$,  the forms
\begin{align*}
\rho_0 &= e^{i\omega}|_{\tilde U_0},\\
\rho_i &= z_i^{a_i} + dz_i\wedge dw_i,\\
g_{i0}(z_i)&=z_i^{a_i},\\
B_{i0}&=B_i|_{\tilde U_i\cap\tilde U_0}.
\end{align*}
define a generalized complex structure on the total space
$\Omega^1_\Sigma(D)$ of the cotangent bundle of $\Sigma$ twisted by
the line bundle with Euler class Poincar\'{e} dual to $D$.

When $\Sigma$ is non-orientable, this generalized complex structure
depends only on the diffeomorphism class of $\Sigma$ and the set
$a=\{a_i,\ldots,a_n\}$ of multiplicities; we denote this generalized
complex manifold by $\Omega^1_\Sigma(a)$, and simply
$\Omega^1_\Sigma(n)$ if all $a_i=1$.

When $\Sigma$ is oriented, however, the chosen complex coordinates
$z_i$ may not be compatible with the orientation; indeed, we may
divide the points into two groups
$\{p_1,\ldots,p_k\},\{p_{k+1},\ldots,p_n\}$ according to whether
$z_i$ is compatible with orientation or not, respectively. Then the
generalized complex structure depends only on the diffeomorphism
class  of the oriented 2-manifold $\Sigma$ and the two sets of
multiplicities $a_+=\{a_1,\ldots, a_k\}$,
$a_-=\{a_{k+1},\ldots,a_n\}$. We denote this generalized complex
manifold by $\Omega^1_\Sigma(a_+,a_-)$, and simply
$\Omega^1_\Sigma(k,n-k)$ if all $a_i=1$.   In this case, therefore,
we obtain a generalized complex structure on $\Omega^1_\Sigma(D)$
such that the zero section $\Sigma$ defines an oriented 2-brane
which intersects the complex loci $\pi^{-1}(p_i),\ i\leq k$
positively and the remaining complex loci $\pi^{-1}(p_i),\ i>k$
negatively, with respect to the natural orientation on the complex
locus.
\end{ex}

We have seen that 2-branes in even generalized complex 4-manifolds
are generically Lagrangian; we now prove that such 2-branes have a
standard tubular neighbourhood up to $B$-diffeomorphism, just as in
the familiar case of Lagrangian submanifolds of symplectic
manifolds.  We restrict to the case where the 2-brane intersects the
complex locus transversally; in this case, the standard
neighbourhood is given by a neighbourhood of the zero section in
$\Omega^1_\Sigma(n)$ or $\Omega^1_\Sigma(k,n-k)$ from
Example~\ref{normalforms}.

\begin{theorem}[Brane \nhood\ theorem]\label{theo:lagrangian nhood}
Let $(\Sigma,F)\stackrel{\iota}{\hookrightarrow}(M,\J)$ be a compact
2-brane in an even generalized complex 4-manifold which intersects
the complex locus transversally at n nondegenerate points.

\begin{itemize}
\item[i)]If $\Sigma$ is non-orientable, then it has a tubular
neighbourhood isomorphic to a \nhood\ of the zero section in
$\Omega^1_\Sigma(n)$.

\item[ii)]Otherwise, orient $\Sigma$ and let $k$ be the number of points
where its intersection with the complex locus is positive. Then
$\Sigma$ has a tubular neighbourhood isomorphic to a \nhood\ of the
zero section in $\Omega^1_\Sigma(k,n-k)$.
\end{itemize}
\end{theorem}

The proof will require the following two preliminary results.

\begin{lem}\label{item:1}
Under the hypotheses of Theorem~\ref{theo:lagrangian nhood}, one can
find coordinates around a point in the intersection of the brane
with the complex locus so that the \gcs\ is determined by
\begin{equation}\label{repform}
\rho = w +dw\wedge dz
\end{equation}
 and the brane is given by $z =0$.
\end{lem}
\begin{proof}
According to Theorem \ref{prop:local form}, a \nhood\ of a
nondegenerate point in the complex locus is $B$-diffeomorphic to a
neighbourhood of the origin in $\C^2$, equipped with the generalized
complex structure given by
\[w+ dw \wedge dz= w\exp(\tfrac{dw\wedge dz}{w}).\]
Choose such coordinates so that the origin is a point $p$ in the
intersection of $\Sigma$ with the complex locus.

Now let $B + i \omega = \tfrac{dw\wedge dz}{w}$.  Since $(\Sigma,F)$
is a brane, $F = B|_{\Sigma}$ is  a well defined 2-form.  Since
$\Sigma$ is transversal to $w=0$ at $p$, in a neighbourhood around
$p$ we can write $F= \iota^*\tilde F$, where $\tilde F = -i
b(w,\overline{w}) dw \wedge d\overline{w}$, for a real function $b$.
Note that $\tilde F$ is a closed 2-form on a \nhood\ of $0$ in
$\C^2$.

Now we can perform a $B$-field transform by $\tilde F$ to obtain:
\begin{align*}
e^{-\tilde F}\rho &= w +dw\wedge dz -i w b dw\wedge
d\overline{w}\\
&= w + dw\wedge (dz -iw b d\overline{w})\\
&= w + dw\wedge d\tilde z,
\end{align*}
for new complex coordinates $(w,\tilde z)$. Since $\Sigma$ is a
brane, $\omega|_{\Sigma} = 0$ and by construction $(B-F)|_{\Sigma}
=0$, hence $dw\wedge d\tilde z$ annihilates $T\Sigma$, which shows
that $\Sigma$ is a complex submanifold of $\C^2$ \wrt\ this new
complex structure. Therefore there is a holomorphic function $f(w)$
such that $\Sigma = \{(w, f(w))\}$ and hence $(w,z):=(w, \tilde z
-f(w))$ are the required coordinates in which the generalized
complex structure has the standard form \eqref{repform} and $\Sigma$
is given by $\{z=0\}$.\end{proof}

\begin{lem}\label{item:2}
Under the hypotheses of Theorem~\ref{theo:lagrangian nhood}, the
normal bundle to $\Sigma$ is diffeomorphic to $\Omega^1_\Sigma(D)$,
i.e. the cotangent bundle twisted by the complex line bundle
Poincar\'{e} dual to $[D]\in H_0(\Sigma,\Z_\omega)$  (homology with
orientation-twisted coefficients).  Here $D=\sum p_i$, for
$\{p_1,\ldots,p_n\}$ the intersection of $\Sigma$ with the complex
locus.
\end{lem}

\begin{proof}
The bundle $\tau_F$ is an extension of the form
\[0 \into N^*\Sigma \into \tau_F \into T\Sigma  \into 0.\]
composing $\J$ with the projection onto the tangent bundle $\pi_T$,
we obtain a map
\[\pi_T \circ \J: \N^* \into T\Sigma.\]
Fibrewise, this map vanishes precisely at the intersection points
$\{p_i\}$ of $\Sigma$ with the complex locus; otherwise it is an
isomorphism.

By Lemma \ref{item:1}, near an intersection point, we can find
coordinates so that the structure is given by \eqref{repform} and
$\Sigma$ is given by $\{z =0\}$.  In this case, $dz$ is a section of
$N^*\Sigma$ and $\pi \J dz = 2i w\del_{w}$. Therefore this point
contributes with a $+1$ to the Euler characteristic of $T\Sigma$,
i.e.
\[\chi(T \Sigma) = \chi(N^*\Sigma) + n.  \]
Hence, as differentiable bundles, we have
\[N\Sigma \cong \Omega^1_\Sigma(D).\]
\end{proof}

\begin{proof}[Proof of Theorem \ref{theo:lagrangian nhood}]

By Lemma~\ref{item:2}, a tubular neighbourhood of $\Sigma$ is
diffeomorphic to a neighbourhood of the zero section in
$\Omega^1_\Sigma(D)$, for $D=\sum p_i$ given by the sum of the
intersection points with the complex locus.  This diffeomorphism can
be chosen so that the generalized complex structure $\J$ agrees, at
the points $\{p_i\}$, with the normal form $\J_0$  given in
Example~\ref{normalforms}, namely $\Omega^1_\Sigma(n)$ (for case
$i)$) or $\Omega^1_\Sigma(k,n-k)$ (for case $ii)$).

According to Lemma~\ref{item:1}, there is a $B$-diffeomorphism
identifying $\J$ with $\J_0$ in neighbourhoods of the intersection
points $\{p_i\}$.  Away from these neighbourhoods, $\J,\J_0$ have
the form $\exp(B+i\omega)$, $\exp(B_0+i\omega_0)$, and  $\Sigma$ is
simply a Lagrangian submanifold for the symplectic structure
$\omega$, by Proposition~\ref{prop:lagrangian brane}. Moser's
argument then furnishes a diffeomorphism $\psi$, compactly supported
in the symplectic locus and fixing $\Sigma$, and such that
$\psi^*\omega = \omega_0$ in a tubular neighbourhood of $\Sigma$.
Finally we apply the $B$-field transform by $B_0 - \psi^*B$ to
identify $\J$ with $\J_0$ on the tubular neighbourhood, as required.
\end{proof}


\section{Blowing up and down}\label{sec:blow-up}

A common feature of complex and symplectic manifolds is that any
point $p$ may be blown up to obtain a new complex or symplectic
manifold, where $p$ has been replaced by a complex projective space
of real codimension 2, called the exceptional divisor. In the
complex case, the blowup is uniquely determined by the choice of
point; the exceptional divisor may be canonically identified with
the projectivized tangent space to $p$.  In the symplectic case, the
blowup is not unique: it depends on a real parameter measuring the
symplectic size of the exceptional divisor.

A point of symplectic type in a generalized complex manifold has a
neighbourhood $B$-diffeomorphic to a symplectic structure.  Hence,
symplectic blow-up may be used to produce new generalized complex
manifolds just as is done in symplectic geometry. Since this
construction is based on the symplectic blow-up, the \gcs\ obtained
is non-unique.

In this section, we show that a nondegenerate complex point $p$ in a
generalized complex 4-manifold $M$ may be blown up in a canonical
fashion, just as for a complex manifold. Since the tangent space
$T_pM$ is complex, we may identify the exceptional divisor $\Sigma$
with the complex projective line $\C\mathbb{P}(T_p M)$, which
contains a distinguished point $\tilde p$ corresponding to the
tangent line to the locus of complex points near $p$. A
neighbourhood of $\Sigma$ is then isomorphic to a neighbourhood of
the zero section in the tautological line bundle over $\Sigma$,
equipped with the generalized complex structure from
Example~\ref{ex:line bundles} with $D = \tilde p$; that is, $\Sigma$
becomes a 2-brane which is Lagrangian away from $\tilde p$.

Conversely, we show that any 2-brane $\Sigma$ in a generalized
complex 4-manifold $\tilde M$ which intersects the complex locus
transversally in a single nondegenerate point $\tilde p$ may be
blown down, yielding a generalized complex 4-manifold $M$, with a
marked nondegenerate complex point $p$.  This is analogous to the
result that any rational $-1$-curve in a complex surface may be
blown down.

By Theorem~\ref{prop:local form}, a nondegenerate complex point
$p\in M$ has a coordinate neighbourhood $U$ with generalized complex
structure $\J_\beta$ given, in complex coordinates $(w,z)$, by
\begin{equation}\label{stdblow}
\rho = w + dw\wedge dz = e^\beta(dw\wedge dz),
\end{equation}
for the bivector $\beta = w \del_{w}\wedge \del_{z}$.  Let $\tilde
U$ be the usual complex blowup of $U$ at $p$, so that $\tilde U$ may
be described as a neighbourhood of the zero section $\Sigma$ of
$\OO(-1)$ over $\C P^1$. The anticanonical section $\beta$ naturally
lifts to a bivector $\tilde\beta$ on the blowup, which then defines
a generalized complex structure $\J_{\tilde\beta}$ on $\tilde U$ as
in Example~\ref{ex:linb}. Furthermore, as in Example~\ref{ex:line
bundles}, the exceptional divisor $\Sigma$ becomes a 2-brane in
$\tilde U$ which is Lagrangian except at its intersection with the
complex point $\tilde p = [z_1=0]$.

Choosing affine charts $\tilde U_1 = \{(w,\tilde z)=(w, z/w)\}$,
$\tilde U_2 = \{(\tilde w,z)=(w/z, z)\}$  for the blowup, the
bivector $\tilde\beta$ may be written as
\[
\tilde\beta = \begin{cases} \del_{w}\wedge \del_{\tilde z} &
\text{in
}\tilde U_1\\
\tilde w\del_{\tilde w}\wedge\del_{z}& \text{in }\tilde U_2,
\end{cases}
\]
giving a generalized complex structure $\tilde \J$ on $\tilde U$
described by the forms
\begin{equation}\label{blrho}
\tilde\rho = \begin{cases} 1 + dw\wedge d\tilde z& \text{in
}\tilde U_1\\
\tilde w + d\tilde w\wedge dz& \text{in }\tilde U_2.
\end{cases}
\end{equation}
This allows us to see explicitly that the exceptional divisor is a
generically Lagrangian 2-brane (since $dw\wedge d\tilde z$ is a
$(2,0)$-form, and hence vanishes upon pullback to $\Sigma$).

\begin{proposition}
Let $\pi:\tilde U\rightarrow U$ be the blowup projection map.  With
respect to the generalized complex structures described above, the
generalized map $\Pi=(\pi,0)$ is holomorphic. Further, its
restriction $\tilde U\backslash{\pi^{-1}(p)}\rightarrow
U\backslash{\{p\}}$ is an isomorphism.
\end{proposition}
\begin{proof}
Let $\pi_i:\tilde U_i\rightarrow U$ be the blowup projection
restricted to the $\tilde U_i$, so that $\pi_1:(w,\tilde z)\mapsto
(w,w\tilde z)$ while $\pi_2:(\tilde w,z)\mapsto (\tilde w z, z)$.
Calculating the pullback of $\rho=w+dw\wedge dz$, we have
\begin{align*}
\pi_1^*\rho &= w + wdw\wedge d\tilde z = w(1+dw\wedge d\tilde z)\\
\pi_2^*\rho &= \tilde w z + zd\tilde w\wedge dz= z(\tilde w + \tilde
w\wedge dz).
\end{align*}
Comparing with \eqref{blrho}, we see that $\pi^*\rho$ defines the
same generalized complex structure as $\tilde \rho$, away from the
exceptional divisor.  Since holomorphicity is a closed condition, we
conclude that the generalized map $(\pi,0)$ is holomorphic, as
required.
\end{proof}

\begin{proposition}\label{canon}
Any generalized complex automorphism $\Phi=(\varphi,B)$ of $U$
fixing $p$ has a canonical lift to an automorphism $\tilde\Phi =
(\tilde\varphi,\tilde B)$ of the blowup $\tilde U$ making the
following diagram commute, where $\Pi$ denotes the generalized
holomorphic projection $(\pi,0)$.
\begin{equation}\label{canonsq}
\xymatrix{\tilde
U\ar[r]^{\Pi}\ar[d]_{\tilde\Phi}&U\ar[d]^{\Phi}\\\tilde
U\ar[r]_{\Pi}&U}
\end{equation}
\end{proposition}
\begin{proof}
Since $\Phi$ is a generalized complex automorphism which fixes $p$,
and since $p$ is a complex point, $d\varphi|_p$ must be a complex
linear automorphism of $T_pU$ (and $B|_p$ must be of type $(1,1)$).
Hence $\varphi$ lifts to a diffeomorphism $\tilde\varphi:\tilde
U\rightarrow \tilde U$ which acts via $\mathbb{P}(d\varphi|_p)$ on
the exceptional divisor $\mathbb{P}(T_pU)$ and which coincides with
$\varphi$ elsewhere. Then we may take $\tilde B = \pi^*B$, yielding
$\Pi\circ \tilde\Phi = \Phi\circ\Pi$, and $\tilde \Phi$ must be an
automorphism since $\Pi^{-1}\circ\Phi\circ\Pi$ is an isomorphism on
the dense set $\pi^{-1}(U\backslash\{p\})$.
%
\end{proof}


\begin{theorem}[Blowing up]
For any nondegenerate complex point $p\in M$ in a generalized
complex 4-manifold, there exists a generalized complex 4-manifold
$\tilde M$ and a generalized holomorphic map $\pi:\tilde
M\rightarrow M$ which is an isomorphism $\tilde
M\backslash\pi^{-1}(p)\rightarrow M\backslash{\{p\}}$ and which is
equivalent to $\pi:\tilde U\rightarrow U$ as above in a
neighbourhood of $\pi^{-1}(p)\cong\C \mathbb{P}(T_pM)$.  The pair
$(\tilde M, \pi)$ is called the blowup of $M$ at $p$, and is unique
up to canonical isomorphism.
\end{theorem}
\begin{proof}
Choose a coordinate neighbourhood $U$ of $p$ which is standard in
the sense of~\eqref{stdblow} and let $\pi:\tilde U\rightarrow U$ be
the standard complex blowup as above.  Then define $\tilde M =
M\backslash\{p\}\cup_{\pi}\tilde U$ and extend $\pi$ by the identity
map to all of $\tilde M$.  Then by Proposition~\ref{canon}, $\tilde
M$ is canonically independent of the chosen coordinates.
\end{proof}

By Theorem~\ref{theo:lagrangian nhood}, any 2-brane which intersects
the complex locus in a single nondegenerate point has a standard
tubular neighbourhood; in particular, such a 2-brane must have
self-intersection $-1$. Observing that the exceptional divisor
$\Sigma = \pi^{-1}(p)$ of the blowup described above is precisely
such a 2-brane, we obtain the following result.
\begin{theorem}[Blowing down]\label{blow-up}
A generalized complex 4-manifold $\tilde M$ containing a 2-brane
$\Sigma\cong S^2$ intersecting the complex locus in a single
nondegenerate point may be blown down; i.e. there is a generalized
holomorphic map $\pi:\tilde M\rightarrow M$ to a generalized complex
manifold $M$ which is an isomorphism $\tilde
M\backslash\Sigma\rightarrow M\backslash\{p = \pi(\Sigma)\}$, and
which is equivalent to $\pi:\tilde U\rightarrow U$ as above in a
neighbourhood of $\Sigma$.
\end{theorem}


\section{$C^\infty$ log transform}\label{sec:surgery}

In \cite{CG06}, the authors introduced a construction of \gcms,
wherein a symplectic $4$-manifold $(M,\omega)$ undergoes surgery
along an embedded symplectic 2-torus with trivial normal bundle to
yield a \gcm\ with type change along a 2-torus.  This type of
4-manifold surgery is called a $C^\infty$ log
transform~\cite{GoMo93}. We now clarify this construction and study
its effect on a Lefschetz fibration.

Let $T\hookrightarrow M $ be a symplectic 2-torus with trivial
normal bundle and symplectic area $A$.  By Moser's argument, we may
choose polar coordinates $(r,\theta_1)$ transverse to $T$ and
angular coordinates $(\theta_2,\theta_3)$ along $T$ such that the
symplectic form becomes, for $r<\eps$,
\begin{equation}\label{sympfo}
\omega = rdr\wedge d\theta_1 + \tfrac{A}{4\pi^2} d\theta_2\wedge
d\theta_3.
\end{equation}
Let $U_0=M\backslash\{(r,\theta_1,\theta_2,\theta_3)\ :\
r\leq\tfrac{\eps}{2}\}$, and equip it with the given symplectic
structure
\[
\rho_0 = e^{i\omega}|_{U_0}.
\]
We now compare this to the singular symplectic structure near the
complex locus of a generalized complex manifold.

Consider the quotient of the generalized complex structure from
Example~\ref{typechang} by the lattice $\Gamma = \langle
A,iA\rangle\subset\C$, where $\gamma\in\Gamma$ acts via
$\gamma:(z,w)\mapsto (z+\gamma,w)$.  Then the complex locus $\Sigma
= \{w=0\}$ is an elliptic curve with modular parameter $\tau=i$. As
per Corollary~\ref{cystr}, $\Sigma$ inherits a canonical holomorphic
differential $\Omega=dz$, which in turn defines $A$ via $A^2 =
\int_\Sigma i\Omega\wedge\bar\Omega$. With respect to polar
coordinates $w = \tilde r e^{i\tilde\theta_1}$ and angular
coordinates $\tilde \theta_2= \tfrac{2\pi}{A}\mathrm{Re}(z), \tilde
\theta_3 = \tfrac{2\pi}{A}\mathrm{Im}(z)$, the generalized complex
structure is given by
\[
\rho_1 = \tilde r e^{i\tilde\theta_1}(1 + \tfrac{A}{4\pi^2}(d\log
\tilde r + id\tilde\theta_1)\wedge(d\tilde \theta_2 + i
d\tilde\theta_3)),
\]
which is proportional to $e^{\tilde B + i\tilde\omega}$ for the
singular forms
\begin{align*}
\tilde B &= \tfrac{A }{4\pi^2}(d\log \tilde r \wedge d\tilde\theta_2 -  d\tilde \theta_1\wedge d\tilde\theta_3),\\
\tilde\omega &= \tfrac{A}{4\pi^2}(d\log \tilde r \wedge d\tilde
\theta_3 + d\tilde \theta_1\wedge d\tilde\theta_2).
\end{align*}
We then observe that $\tilde\omega$ (for $\tilde r>\tilde r_0>0$) is
symplectomorphic to $\omega$ (for $r>0$) via the coordinate
transformation
\[
\varphi_{\tilde r_0}: (r,\theta_1,\theta_2,\theta_3)\mapsto (\tilde
r,\tilde \theta_1,\tilde \theta_2,\tilde\theta_3)= ({\tilde r_0}
e^{{2\pi^2 A^{-1}}r^2},\theta_2,\theta_3, \theta_1).
\]
If we glue the open set $U_0$ to the open set $U_1 = \{(\tilde
r,\tilde\theta_1,\tilde\theta_2,\tilde\theta_3)\ :\ \tilde r < 1\}$
along the neck $\eps/2<r<\eps$, using the diffeomorphism
$\varphi_{\tilde r_0}$, with ${\tilde
r_0}=e^{-2\pi^2A^{-1}{\eps}^2}$ (so that, for $\eps/2<r<\eps$ we
have ${\tilde r_0} e^{\tfrac{\pi^2\eps^2}{2A}}<\tilde r<1$), then we
observe that, in $U_0\cap U_1$, we have
\[
e^{B_{01}}\varphi_{\tilde r_0}^*\rho_1 = e^{i\omega} = \rho_0,
\]
for the closed 2-form $B_{01}\in \Omega^{2}(U\cap V,\R)$ given by
\begin{equation}\label{cechh}
B_{01}=-rdr\wedge d\theta_3 - \tfrac{A}{4\pi^2}d\theta_1\wedge
d\theta_2.
\end{equation}
Choosing a \v{C}ech trivialization $B_{01} = (B_1 - B_0)|_{U_0\cap
U_1}$ for $B_i\in \Omega^2(U_i,\R)$, we see that $\{e^{B_i}\rho_i\}$
defines a generalized complex structure on the surgered manifold
$\hat M = U_0\cup_{\varphi_{\tilde r_0}} U_1$, which has symplectic
type in $U_0$ and which changes type in $U_1$ along an elliptic
curve. While the underlying 3-form of the original symplectic
manifold $(M,\omega)$ vanishes, this is not the case for $\tilde M$,
where the 3-form is given by $H|_{U_i} = -dB_i$.

The cohomology class $[H]$ may be easily described, since for any
closed 1-form $\xi$, we have
\[
\int_{\hat M} H\wedge\xi = -\int_{T^3\simeq U_0\cap U_1}
B_{01}\wedge\xi.
\]
Hence by~\eqref{cechh}, $[H]$ is Poincar\'e dual to $A$ times the
circle parametrized by $\theta_3$. Since $\theta_3$ corresponds to
$\tilde \theta_2$ in the gluing, we may describe this circle as an
integral circle of the real part of the canonical holomorphic vector
field $Z = \del/\del z$ on $\Sigma$. The 4-manifold surgery
described above is known as a $C^\infty$ logarithmic transform of
multiplicity zero. We now summarize the above discussion.

\begin{theorem}[Cavalcanti--Gualtieri \cite{CG06}]\label{theo:surgery}
Let $(M,\omega)$ be a symplectic 4-manifold, $T \hookrightarrow M$
be a symplectic 2-torus of area $A$ with trivial normal bundle. Then
the multiplicity zero $C^{\infty}$ logarithmic transform of $M$
along $T$, denoted $\hat M$,  admits a \gcs\ such that:

\begin{enumerate}
\item  The complex locus is given by an elliptic curve $\Sigma$ with modular
parameter $\tau=i$, and the induced holomorphic differential
$\Omega$ has periods $\langle A, iA\rangle$.

\item Integrability holds with respect to
a 3-form $H$, which is Poincar\'e dual to $A$ times the homology
class of an integral circle of $\Re(\Omega^{-1})$ in $\Sigma$.
\end{enumerate}
\end{theorem}

\begin{ex}[Lefschetz fibrations]
Let $(M,\omega)$ be a symplectic 4-manifold which is expressed as a
symplectic Lefschetz fibration (in the sense of~\cite{Don}) $\pi:
M\rightarrow B$ whose generic fiber has genus 1. If $T$ is a smooth
fiber, then it is a symplectic 2-torus with trivial normal bundle.
Trivialize the fibration near $T$ so that a \nhood\ of $T$ is
diffeomorphic to $D^2 \times T^2$, where $\pi$ is the first
projection. In coordinates, we have
$\pi(r,\theta_1,\theta_2,\theta_3) = (r,\theta_1)$.

We may now apply an isotopy, supported in a sufficiently small
neighbourhood of $T$, taking $\omega$ into the standard
form~\eqref{sympfo} in some neighbourhood of $T$ and such that $\pi$
remains a symplectic Lefschetz fibration.

Now, let $\hat{M}$ be the $C^{\infty}$ logarithmic transform of $M$
along $T$ as in Theorem \ref{theo:surgery}. Then, using the
preceding notation, the fibration projection $\pi$ is still well
defined on $U_0$.  We extend $\pi$ to $U_1$, and hence to all of
$\hat{M}$, as follows.  For points in $U_0\cap U_1$ (for which
$\tilde r> {\tilde r_0} e^{\pi^2\eps^2(2A)^{-1}}$), the projection
may be written

\begin{equation}\label{rtr}
\pi:(\tilde r, \tilde\theta_1,\tilde\theta_2,\tilde\theta_3)\mapsto
(r,\theta_1)=((\tfrac{A}{2\pi^2}\log \tfrac{\tilde r}{{\tilde
r}_0})^{1/2},\tilde\theta_3).
\end{equation}
To extend this to all of $U_1$, we choose $f:[0,1]\rightarrow
[\eps/4,\eps]$ to be a smooth monotone function of $\tilde r$ such
that $f(0)=\eps/4$ and $f$ agrees with $r = r(\tilde r)$
(see~\eqref{rtr}) for $ {\tilde r_0} e^{(2A)^{-1}\pi^2\eps^2}<\tilde
r<1$.  Then we define, for all $\tilde r\in [0,1]$,
\[
\hat\pi|_{U_1}:(\tilde
r,\tilde\theta_1,\tilde\theta_2,\tilde\theta_3)\mapsto (f(\tilde
r),\tilde\theta_3).
\]
This defines a projection of $U_1$ onto an open-closed annulus
$\eps/4\leq r<1$; together with the projection $\pi:U_0\rightarrow
B\backslash\{r\leq \eps\}$, we obtain a projection of $\hat M$ onto
the surface with boundary $\hat B = B\backslash D^2_{\eps/4}$:
\[
\hat\pi:\hat M\longrightarrow \hat B = B\backslash\{(r,\theta)\ :\
r<\eps/4\}.
\]
This map is a symplectic Lefschetz fibration away from the boundary
$\del \hat B = \{r=\eps/4\}$, where the fiber degenerates to a
circle; indeed $\hat \pi^{-1}(\del \hat B)$ is precisely the
elliptic curve forming the complex locus of the generalized complex
manifold.  The fibers over boundary points are the integral circles
for $\Re(Z)$, and we obtain from $\hat\pi_*\Im(Z)$ a vector field
along the boundary with period $A$.  Note also that since we are not
changing the Lefschetz fibration outside a neighbourhood of
$\del\hat B$, the monodromy about any path homotopic to the boundary
is trivial.

The generalized Lefschetz fibration described above provides a
pictorial description of the behaviour of the generalized complex
structure: the geometry is symplectic over the interior $\hat
B\backslash \del B$ and complex over the boundary $\del B$.
Furthermore, the generalized complex structure is integrable \wrt\ a
closed 3-form $H$ which is Poincar\'e dual to $A$ times the homology
class of the circle $\hat\pi^{-1}(p)$, for $p\in \del \hat B$.

\begin{figure}[h!!]
\begin{center}
{\psfig{file=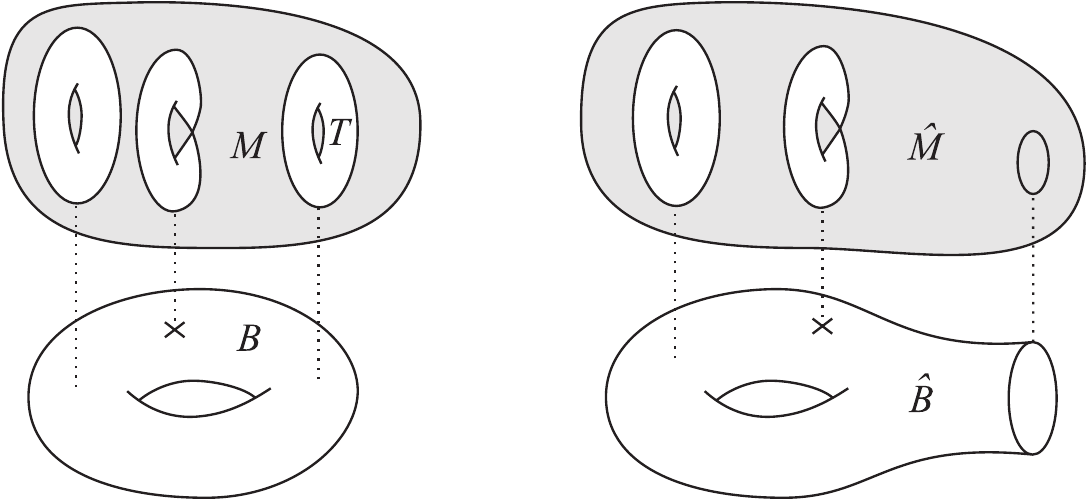, height=3.5cm,clip=}} \caption{A symplectic
genus 1 Lefschetz fibration undergoes surgery along a smooth fiber
$T$, becoming a generalized Lefschetz fibration of a generalized
complex 4-manifold over a surface with boundary.} \label{sboundary}
\end{center}
\end{figure}
\end{ex}

\section{Examples}\label{sec:examples}
In this section, we use the tools introduced above, namely blowing
up and down as well as the $C^\infty$ log transform, in conjunction
with the representation of a generalized complex 4-manifold as a
generalized elliptic Lefschetz fibration, to produce new examples of
generalized complex manifolds.  In particular, we show that the
connected sum of any odd number of copies of $\C P^2$ has a \gcs.

\begin{ex}[Fiber sums]\label{fiber}
Given symplectic manifolds $M_1, M_2$, each equipped with a genus 1
Lefschetz fibration over bases $B_1,B_2$ respectively, we may
produce their symplectic fiber sum (see ~\cite{Gom}), denoted
$M_1\#_{f}M_2$, using a symplectic identification of smooth fibers.
The fiber sum is then a Lefschetz fibration over the connected sum
$B_1\# B_2$.

The $C^{\infty}$ log transform $\widehat{M_1 \#_{f} M_2}$ then has a
generalized Lefschetz fibration over a manifold with boundary which
is precisely the \emph{boundary connected sum} of the surfaces with
boundary $\hat B_1,\hat B_2$ which form the bases of the generalized
Lefschetz fibrations associated to the $C^\infty$ log transforms
$\hat M_1, \hat M_2$. We therefore obtain a connected sum operation
for the generalized Lefschetz fibrations described above along the
boundary fiber.

\begin{figure}[h!!]
\begin{center}
{\psfig{file=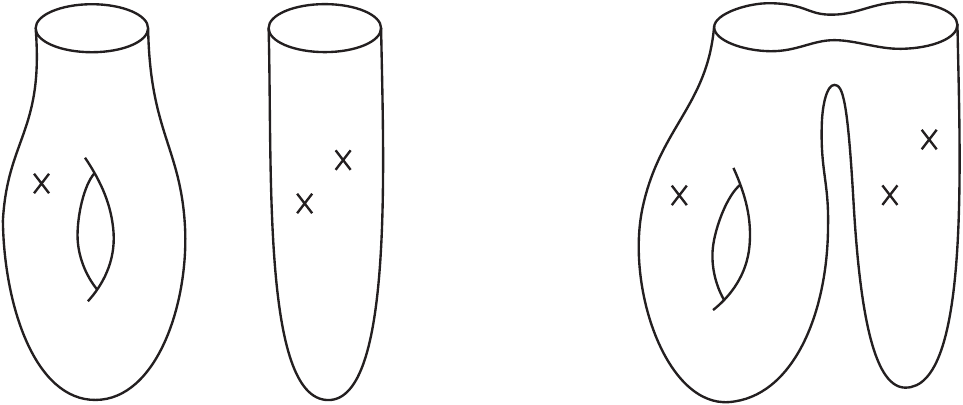, height=2.8 cm,clip=}} \caption{The
connected sum of generalized Lefschetz fibrations over $\hat
B_1,\hat B_2$ along the boundary fiber is itself a generalized
Lefschetz fibration over the boundary connected sum $\hat B_1\#_\del
\hat B_2$.}
\end{center}
\end{figure}


\end{ex}

\begin{ex}[Branes and blow down]\label{ex:blow down}
A standard tool in symplectic topology, described in~\cite{Don}, is
the construction of Lagrangian spheres in symplectic manifolds by
the method of vanishing cycles.  This proceeds essentially by
choosing a path in the base of a Lefschetz fibration connecting two
nodal fibers which is such that the same cycle degenerates at each
end of the path.  Using a connection determined by the symplectic
orthogonal of the fibers, a representative of the vanishing cycle
traces out a Lagrangian sphere fibering over the original path.

For a generalized complex 4-manifold presented as a generalized
Lefschetz fibration, we may use the method of vanishing cycles to
construct two types of branes in addition to the Lagrangian spheres
which exist in the symplectic locus.  The first is obtained by
connecting two boundary points in the base by a path such that the
same cycle degenerates near each end.  In this case we obtain a
2-sphere which is Lagrangian in the symplectic locus and intersects
the complex locus transversally at two points. Proposition
\ref{prop:lagrangian brane} then implies that this sphere is a
brane, which according  to Theorem \ref{theo:lagrangian nhood} has
trivial normal bundle.

The more interesting case arises when connecting a boundary point in
the base with the base point of a nodal fiber in such a way that the
same cycle degenerates at each end.  Then we obtain a 2-sphere brane
 intersecting the complex locus transversally at precisely one
nondegenerate point (hence by Theorem \ref{theo:lagrangian nhood} it
has self-intersection $-1$). According to Theorem \ref{blow-up} we
may blow down such a spherical brane to obtain a new \gcm .  Note
that the blow down does not inherit a Lefschetz fibration; rather,
it may be viewed as an analogue of a Lefschetz pencil.

\begin{figure}[h!!]
\begin{center}
{\psfig{file=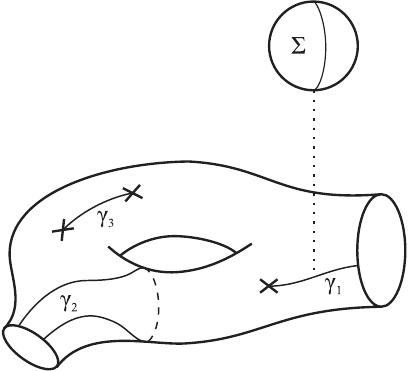, height=4 cm,clip=}} \caption{The path
$\gamma_1$ joins a boundary fiber to a nodal fiber with equal
vanishing cycle
 and lifts to a 2-sphere brane $\Sigma$ with self-intersection $-1$.
The path $\gamma_2$ joins boundary fibers with equal vanishing cycle
and lifts to a 2-sphere brane with trivial normal bundle.  The path
$\gamma_3$ joins nodal fibers with same vanishing cycle, and lifts
to a usual Lagrangian sphere in the symplectic locus.}
\end{center}
\end{figure}

\end{ex}

In the final example, we show that the connected sum of any odd
number of copies of $\C P^2$ admits a \gcs.   This is of particular
interest for the following reasons. First, due to the Kodaira
classification of complex surfaces and a fundamental result in
Seiberg--Witten theory, the manifolds $n \C P^2 \#m\overline{\C
P^2}$ have neither complex nor symplectic structures when $n>1$
\cite{Kod64,Ta94,Wi94}. Therefore, this example produces a family of
\gcms\ which do not admit complex or symplectic structures. Second,
due to results of Hirzebruch and Hopf \cite{HH58}, an oriented
simply connected 4-manifold admits an almost complex structure if
and only if $b_2^+$ is odd.
Since generalized complex manifolds are necessarily almost
complex~\cite{Gu7}, this example shows that for connected sums of
$\C P ^2$ and $\overline{\C P^2}$, there is no obstruction to the
existence of a \gcs\ besides that of being almost complex.

\begin{ex}[A \gcs\ on $(2n-1)\C P^2$]\label{ex:2n-1CP2}
The blow up of $\C P^2$ at the 9 points of intersection of two
generic cubics provides a basic example of a symplectic $4$-manifold
equipped with an elliptic Lefschetz fibration.  This manifold is
sometimes denoted by $E(1)$, and $E(n)$ is the fiber sum of $n$
copies of $E(1)$.  For example, $E(2)$ is diffeomorphic to a $K3$
surface.  A relevant fact concerning the manifolds $E(n)$ is that
the multiplicity zero $C^{\infty}$ logarithmic transform along a
smooth fiber, $\widehat{E(n)}$, ``dissolves'' \cite{Go99}:
\[\widehat{E(n)} = (2n-1) \C P^2 \# (10n -1)\overline{\C P^2}.\]
For $n=1$ we obtain the diffeomorphism $\widehat{E(1)} = \C P^2\# 9
\overline{\C P^2}= E(1)$.  In light of Example \ref{ex:blow down},
we can see the nine $-1$ spheres in the generalized complex manifold
$\widehat{E(1)}$ explicitly in the following way.  The Lefschetz
fibration of $E(1)$ over $S^2$ has twelve nodal fibers with critical
values which we label $\{y_1,y_2,y_3,x_1,\ldots,x_9\}$.  Choose a
smooth fiber $F_0$ about which to perform the $C^\infty$ log
transform.  With respect to a fixed set of paths joining $F_0$ to
the nodal fibers, and using a basis $\{a,b\}$ for $H_1(F_0,\Z)$, the
degenerating cycles are  $a+3b, a,a-3b$ for the fibers over
$y_1,y_2,y_3$, and $b$ for the remaining 9 fibers over $x_i$ (see
e.g., Example 8.2.11 in \cite{Go99}).  We now perform the $C^\infty$
log transform about $F_0$, collapsing the cycle $b$.  We obtain a
generalized Lefschetz fibration of $\widehat{E(1)}$ over a disk with
boundary, as in Figure~\ref{e1}.

\begin{figure}[h!!]
\begin{center}
{\psfig{file=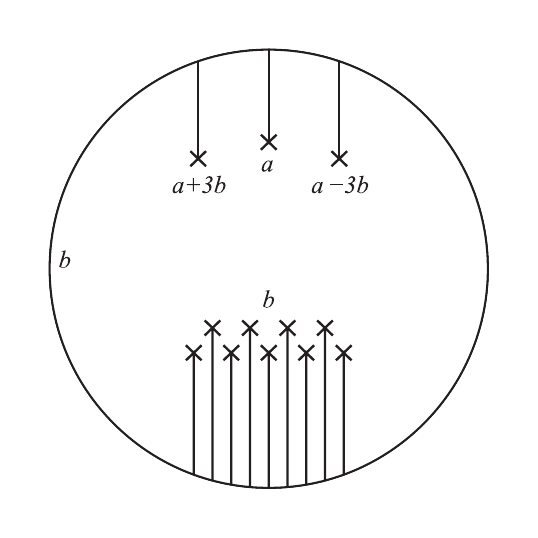, height=5 cm,clip=}} \caption{A
generalized Lefschetz fibration of $\widehat{E(1)}$ over a disk with
boundary, where the $C^\infty$ log transform along $F_0\subset E(1)$
collapses the $b$ cycle at the boundary fiber. The chosen paths are
shown, and the vanishing cycles at nodal fibers are
labeled.}\label{e1}
\end{center}
\end{figure}


Since the paths joining the points $x_i$ to $F_0$ have vanishing
cycle $b$, after the surgery they become paths joining nodal fibers
to boundary points with equal vanishing cycles, and hence lift to a
configuration of nine 2-sphere branes intersecting the complex locus
at single points.  By Theorem \ref{blow-up}, we may blow down each
of these spheres and obtain a \gcs\ on the blow down, which by
Proposition~\ref{kircal} in the Appendix, is diffeomorphic to $\C
P^2$. Generalized complex structures on $\C P^2$ similar to this one
may alternatively be obtained by deforming the complex structure by
a Poisson bivector as in Theorem~\ref{theo:deformations}.

For a more interesting example, construct the connected sum of
$\widehat{E(1)}$ with itself along the boundary fiber, as in
Example~\ref{fiber}, obtaining a new generalized  Lefschetz
fibration over the boundary connected sum, as in
Figure~\ref{cnnect}.
\begin{figure}[h!!]
\begin{center}
{\psfig{file=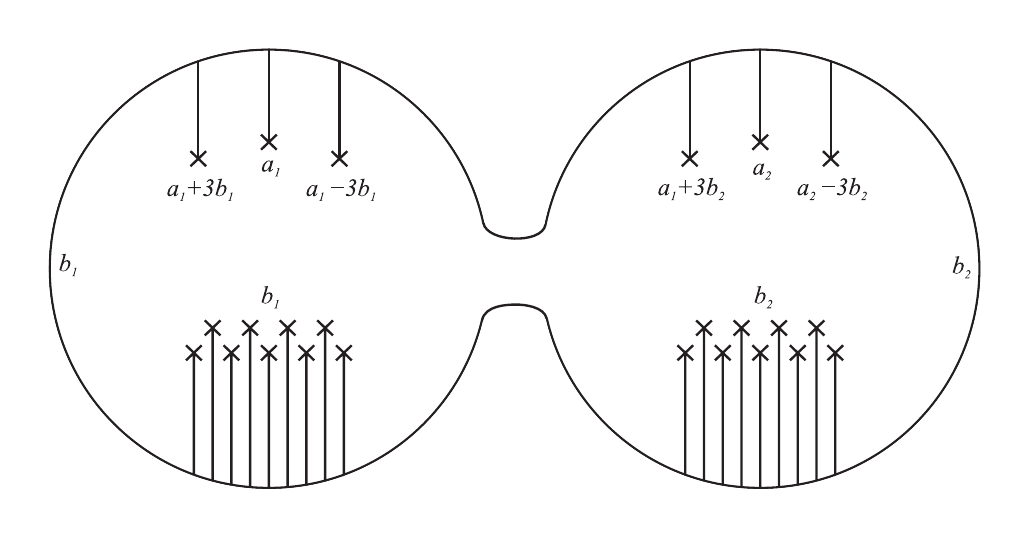, height=5 cm,clip=}}\vspace{-2mm}
\caption{A generalized Lefschetz fibration of $\widehat{E(2)} = 3\C
P^2 \# 19\overline{\C P^2}$, over the boundary connect sum of closed
discs from Figure~\ref{e1}.  The boundary vanishing cycle $b_1$ is
identified with $b_2$ in the connect sum, hence we immediately see
18 2-sphere branes of self-intersection $-1$.}\label{cnnect}
\end{center}
\end{figure}
When performing the fiber sum, we must choose an identification of
the tori which fiber over the boundary circles.  If $(a_i,b_i),
i=1,2$ form a basis for the first homology of a smooth fiber on each
summand, then we require that $b_2$ is identified with $b_1$. We
obtain in this way a generalized complex 4-manifold containing nine
2-sphere branes from each summand, resulting in 18 clearly visible
2-sphere branes which intersect the complex locus in single points.
%
%
There is some freedom in the identification of $a_1$ in the
connected sum; if we set $a_2 = a_1 -7b_1$, then using the monodromy
around $a_1-3b_1$ we can change the cycle $a_2+3b_2 = a_1- 4b_1$
into $-b_1$. The new path which realizes this vanishing cycle is
denoted by $\gamma$ in Figure~\ref{gamma}.
\begin{figure}[h!!]
\begin{center}
{\psfig{file=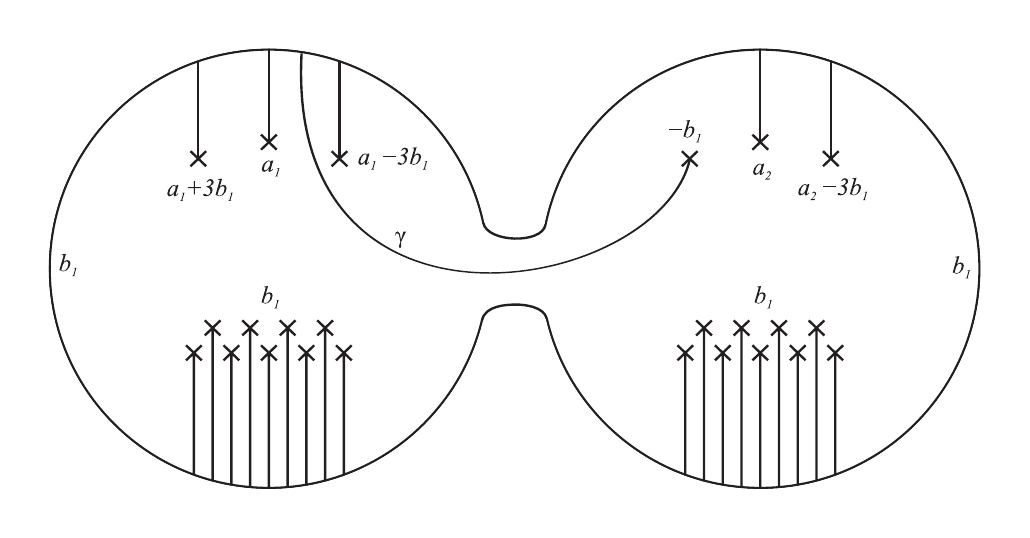, height=5 cm,clip=}}\vspace{-2mm}
\caption{Using the monodromy around the $a_1-3b_1$ critical value,
we reveal a $19^{th}$ 2-sphere brane of self-intersection $-1$,
fibering over the path $\gamma$.}\label{gamma}
\end{center}
\end{figure}

Therefore, in the boundary fiber connected sum $\widehat{E(2)}$, we
find 19 2-sphere branes intersecting the complex locus in single
nondegenerate points.  Since $\widehat{E(2)} = 3\C P^2 \#
19\overline{\C P^2}$, we would like to conclude that upon blowing
down these 19 spheres, we obtain a \gcs\ on the differentiable
manifold $3\C P^2$.  In Proposition~\ref{kircal} (see Appendix), we
use Kirby calculus to verify this claim.

The procedure above may be iterated, taking successive boundary
connect sum with $\widehat{E(1)}$, so that with each new summand we
obtain 10 more 2-sphere branes intersecting the complex locus in
nondegenerate points.  Therefore we obtain $10n -1$ spherical branes
of self-intersection $-1$ in the generalized complex manifold
$\widehat{E(n)}$ which can be blown down.  Using Kirby calculus
again, we prove that the resulting manifold is precisely $(2n-1)\C
P^2$, which therefore has a \gcs.
\begin{figure}[h!!]
\begin{center}
{\psfig{file=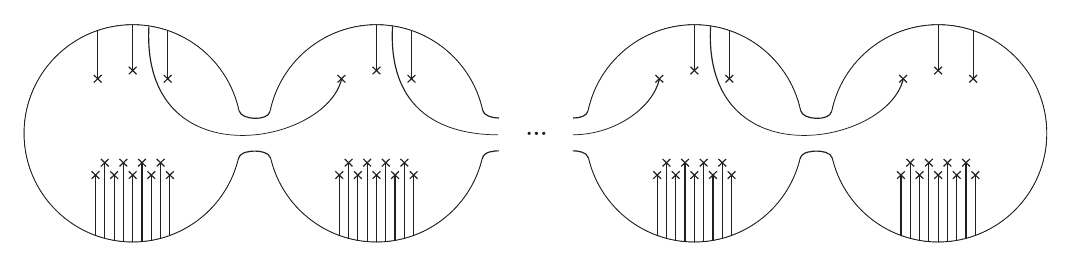, height=3 cm,clip=}}\vspace{-2mm}
\caption{Iterating the boundary fiber connected sum, we find
$(10n-1)$ 2-sphere branes in $\widehat{E(n)}=(2n-1) \C P^2 \# (10n
-1)\overline{\C P^2}$ intersecting the complex locus in single
points. Blowing down each of these, we obtain a \gcs\ on $(2n-1)\C
P^2$.}
\end{center}
\end{figure}
\end{ex}
\begin{remark}
While we use Kirby calculus in the above example to prove that the
manifold obtained after blowing down the $(10n-1)\overline{\C P^2}$
is diffeomorphic to $(2n-1)\C P^2$, there is a simple argument
establishing a homotopy equivalence. Indeed, after blowing down all
the $\overline{\C P^2}$, we obtain a smooth 1-connected manifold
with positive intersection form. By results of Donaldson
\cite{DK98}, the only possible intersection form is the diagonal
$(2n-1) \mathbf{1}$. Since this manifold has the same intersection
form as $(2n-1) \C P^2$, they must be homotopic \cite{Whi49,Mil58}.
\end{remark}

\section*{Appendix 1 -- Kirby calculus}

In the previous section we obtained a manifold, $\widehat{E(n)}$, as
a boundary fiber connected sum of $n$ copies of $\widehat{E(1)}$, or
equivalently, as a $C^\infty$ log transform applied to the
symplectic fiber sum of $n$ copies of $E(1)=\C P^2 \# 9 \overline{\C
P^2}$. In this section we give a handlebody decomposition of
$\widehat{E(n)}$ which makes clear the effect of blowing down the
$(10n-1) \overline{\C P^2}$. This will show that the resulting
manifold is indeed $(2n-1)\C P^2$.

\begin{proposition}\label{kircal}
The manifold  obtained in Example \ref{ex:2n-1CP2} as the blow-down
of the specified $(10n-1) \overline{\C P^2}$ in $\widehat{E(n)}$ is
diffeomorphic to $(2n-1)\C P^2$.
\end{proposition}

\begin{proof}
The first step is to simplify the vanishing cycle information in the
description of $\widehat{E(n)}$ in Figure~\ref{gamma}. Using the
monodromy, we can bundle the $(10n-1)$ singular fibers with
degenerating cycle $b$ together. Furthermore, we can make some
symmetry evident by replacing $a$ by $a+4(n-1)b$. We then obtain the
following vanishing cycles
\begin{equation*}
\begin{cases}
a + (4n-1)b,\\
a + 4kb,\ \text{ for } \ k=n-1,\cdots, 1-n,\\
a-(4n-1)b,\\
 (10n-1)\text{ cycles of type } b,
\end{cases}
\end{equation*}
arranged as in Figure~\ref{all}.
\begin{figure}[h!!]
\begin{center}
{\psfig{file=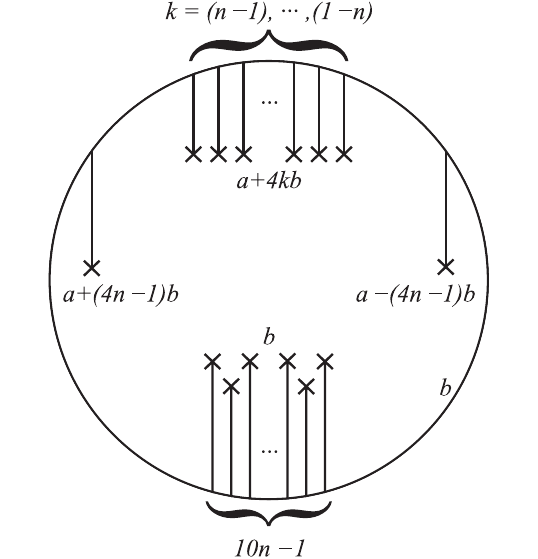, height=5 cm,clip=}}\vspace{-1mm}
\caption{ Generalized Lefschetz fibration of $\widehat{E(n)}$ over a
closed disk.}\label{all}
\end{center}
\end{figure}

The Kirby diagram for $E(n)$ with a regular fiber removed is
obtained by attaching $-1$ framed 2-handles to $T^2 \times D^2$ for
each of the nodal fibers. Then $\widehat{E(n)}$, the surgered
manifold, is obtained from that diagram by attaching a 0-framed
2-handle corresponding to the cycle in $T^2$ which is collapsed
by the surgery \cite{Go99}.  This diagram is shown in
Figure~\ref{fig:kirby1}, where the 2-handles corresponding to the vanishing cycles $a+(4n-1)b, a+4(n-1)b, \cdots, a - 4(n-1)b, a - (4n -1) b$ are denoted by $\alpha_n, \alpha_{n-1},\cdots, \alpha_{-n+1}, \alpha_{-n}$, the 2-handles corresponding to the vanishing cycles $b$ are denoted by $\beta_1, \cdots \beta_{10n-1}$ and the 0-framed 2-handle introduced by the surgery is denoted by $\gamma$ 

\begin{figure}[h!]
\begin{center}
{\psfig{file=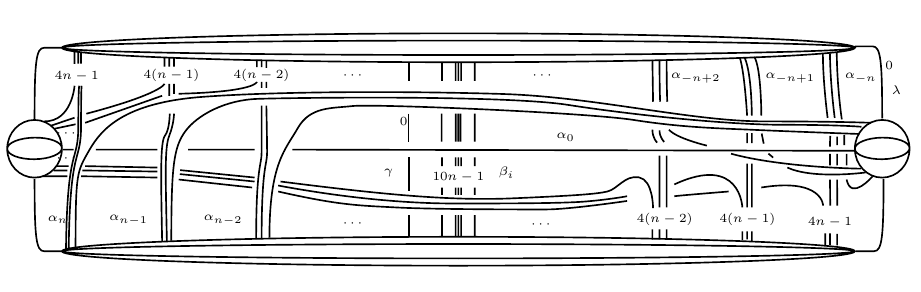, height = 3.22cm, clip=}} \caption{1 and 2-handles of the Kirby diagram for $\widehat{E(n)}$.
All 2-handles  are blackboard $-1$-framed, except for the outermost, $\lambda$,
 and one of the 2-handles running over the vertical 1-handle, $\gamma$. The latter represents the $b$-cycle on the
fibers which collapses on the boundary of the
base.}\label{fig:kirby1}
\end{center}
\end{figure}

If we slide all the handles represented by $\beta_i$ (the circles through the vertically symmetric squashed spheres in Figure
\ref{fig:kirby1}) over the 0-framed 2-handle $\gamma$, also representing the
$b$-cycle, we obtain $(10n-1)$ $-1$-framed unknots which split out
of the diagram. These are precisely the $(10n -1)\overline{\C P^2}$
described in Example~\ref{ex:2n-1CP2}, hence removing them from the
diagram corresponds to blowing down those cycles.  We can also slide
each of the other $-1$-framed 2-handles over a number of copies of $\gamma$ so that the only
2-handle intersecting the top and bottom spheres is $\gamma$ itself. Finally, we can cancel the 1-handle representing $b$
with $\gamma$. If we let $A(n)$ be the manifold obtained after
blowing down, the argument above shows its Kirby diagram is as shown
in Figure~\ref{fig:kirby2}.
\begin{figure}[h!]
\begin{center}
{\psfig{file=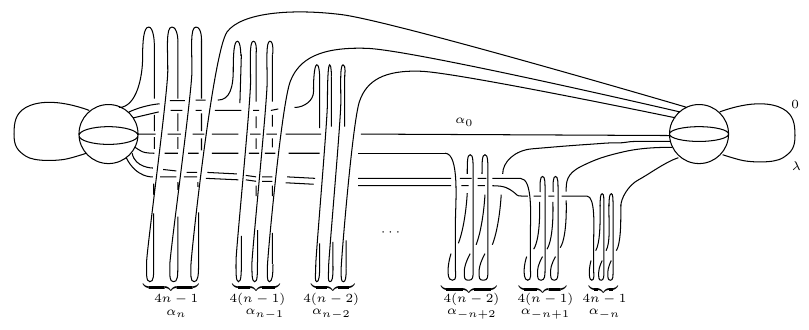, clip=}}
\caption{ 1 and 2-handles of the Kirby diagram of   $A(n)$. All
2-handles  are blackboard $-1$-framed, except for $\lambda$, which is 0-framed.} \label{fig:kirby2}
\end{center}
\end{figure}

Now observe that the 2-handle represented by $\lambda$ can be pushed through the 1-handle and becomes a zero framed unknot disjoint from the rest of the diagram which can therefore be cancelled against a 3-handle.

Then we can slide all the remaining handles through the handle labeled $\alpha_0$ in Figure \ref{fig:kirby2}, so that we can cancel the remaining 1-handle with $\alpha_0$ and  obtain Figure \ref{fig:kirby3} as the Kirby diagram for $A(n)$.
\begin{figure}[h!]
\begin{center}
{\psfig{file=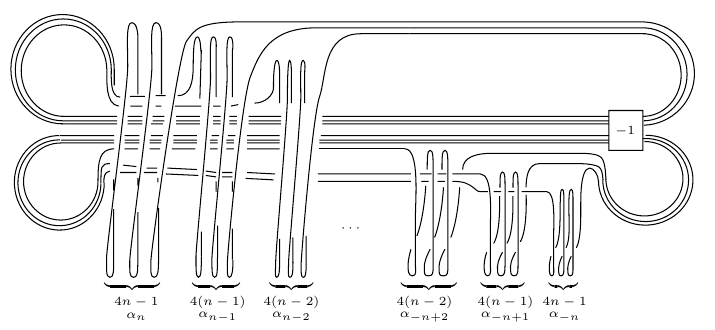, clip=}}
\caption{ $A(n)$ as a 2-handlebody. All 2-handles  are blackboard
$-1$-framed.}\label{fig:kirby3}
\end{center}
\end{figure}

Now that all the 1-handles are gone, we can abandon the blackboard framing and use instead the Seifert framing. Further, we observe that we can move $\alpha_n$ and $\alpha_{-n}$ so that their crossings with the other handles in the diagram are simplified, obtaining Figure \ref{fig:kirby4} as the diagram for $A(n)$.
\begin{figure}[h!]
\begin{center}
{\psfig{file=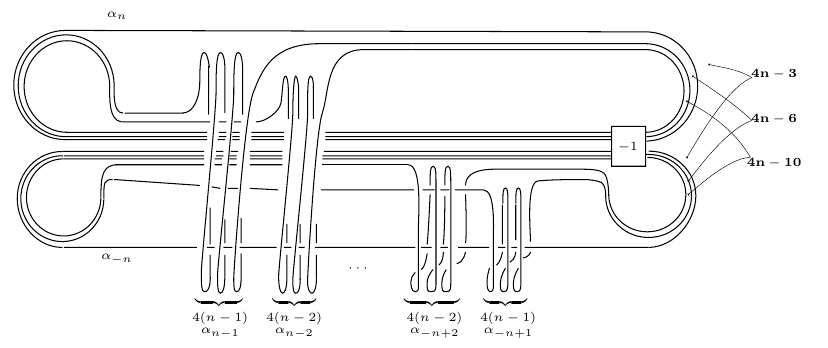, clip=}}
\caption{ $A(n)$ as a 2-handlebody. The boldfaced numbers indicate the canonical framing of the respective handles.}\label{fig:kirby4}
\end{center}
\end{figure}

And then we can move $\alpha_{n-1}$ so that it wraps  around $\alpha_n$ and after that move $\alpha_{n-2}$ so that it wraps around $\alpha_n$ and $\alpha_{n-1}$ and so on. This way, the handles $\alpha_i$ with $i>0$ will only knot with the $\alpha_i$ with $i<0$  in the $-1$ box and we obtain the more symmetric Kirby diagram for $A(n)$ shown in Figure \ref{fig:kirby5}.
\begin{figure}[h!]
\begin{center}
{\psfig{file=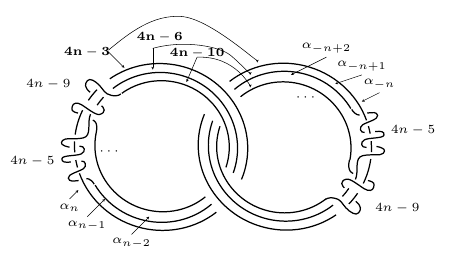, clip=}}
\caption{ $A(n)$ as a 2-handlebody. The boldfaced numbers indicate the canonical framing of the respective handles.}\label{fig:kirby5}
\end{center}
\end{figure}

We now prove by induction that this manifold is $(2n-1)\C P^2$. In
the case when $n=1$, the diagram in Figure~\ref{fig:kirby5} becomes
a pair of 1-framed 2-handles linked in a Hopf link. A handle slide
separates them (see Figure~\ref{fig:kirby6}), rendering a 1-framed
2-handle and a 0-framed 2-handle. The 0-framed 2-handle cancels with
a 3-handle and the result is $\C P^2$, as required.
\begin{figure}[h!]
\begin{center}
{\psfig{file=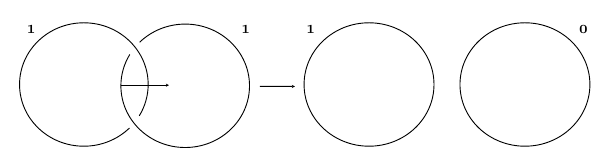, clip=}}
\caption{ A simple handle slide shows that $A(1) = \C
P^2$.}\label{fig:kirby6}
\end{center}
\end{figure}

In the general case, we can slide $\alpha_n$ over $\alpha_{n-1}$ as well as slide $\alpha_{-n}$ over $\alpha_{-n+1}$, turning
them into 1-framed 2-handles knotting only $\alpha_{\pm(n-1)}$. Splitting
out the 1-framed 2-handles (see Figure~\ref{fig:kirby7}), we obtain
the diagram for $2\C P^2 \# A(n-1)$, proving the induction step.
\begin{figure}[h!]
\begin{center}
{\psfig{file=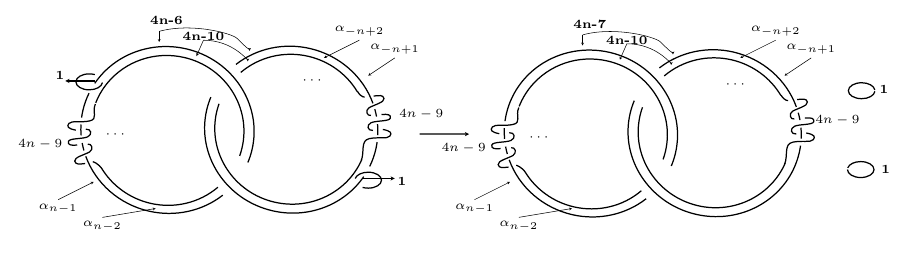, clip=}}
\caption{  In the general case, two handle slides imply $A(n) = A(n-1)
\# 2\C P^2$.} \label{fig:kirby7}
\end{center}
\end{figure}

\end{proof}

\noindent
{\bf Acknowledgements.} We would like to thank Denis Auroux, Nigel Hitchin, David
Martinez and Vicente Mu\~noz for useful conversations.

\end{document}